\documentclass[10pt]{amsart}
\usepackage{amsmath}
\usepackage{amssymb}
\usepackage{amsfonts}
\usepackage{amsthm}
\usepackage{enumerate}
\usepackage[all]{xy}
\usepackage[latin1]{inputenc}
\usepackage{graphicx}
\usepackage{latexsym}
\input xy

\makeatletter
    \newtheoremstyle{definition}
        {5pt}
        {3pt}
        {}
        {0pt}
        {\scshape}
        {.}
        {5pt}
        {\thmname{#1} \thmnumber{#2} \thmnote{[#3]}} 

\newtheoremstyle{theorems}
        {5pt}
        {3pt}
        {\itshape}
        {0pt}
        {\scshape}
        {.}
        {5pt}
        {\thmname{#1} \thmnumber{#2}\thmnote{[#3]}} 

\swapnumbers \theoremstyle{theorems}

\newtheorem{Theo}{Theorem}[section]
\newtheorem{Prop}[Theo]{Proposition}
\newtheorem{PropHS}[Theo]{Proposition (Hattori-Stallings)}
\newtheorem{Cor}[Theo]{Corollary}
\newtheorem{Lemma}[Theo]{Lemma}
\newtheorem{Conj}[Theo]{Extension Conjecture}
\theoremstyle{definition}

\newcommand{\Ext}{{\rm Ext}}

\newcommand{\mmod}{{\rm mod}}

\newcommand{\tr}{{\rm tr}}
\newcommand{\va}{{\varphi}}
\newcommand{\HH}{{\rm H\hspace{-1pt}H}_0}
\newcommand{\La}{\Lambda}

\begin{document}

\title{\sc A proof of the strong no loop conjecture}

\author{Kiyoshi Igusa, Shiping Liu, and Charles Paquette}
\address{Kiyoshi Igusa, Department of Mathematics, Brandeis University, Waltham, MA 02453, United States}
\email{igusa@brandeis.edu}
\address{Shiping Liu, Département de Mathématiques, Université de Sherbrooke, Sherbrooke, Québec, Canada, J1K 2R1}
\email{shiping.liu@usherbrooke.ca}
\address{Charles Paquette, Département de Mathématiques, Université de Sherbrooke, Sherbrooke, Québec, Canada, J1K 2R1}
\email{charles.paquette@usherbrooke.ca}

\begin{abstract}

\vspace{2pt}

The strong no loop conjecture
states that a simple module of finite projective
dimension over an artin algebra has no non-zero
self-extension. The main result of this paper establishes this
well known conjecture for finite dimensional algebras over an
algebraically closed field.

\vspace{-20pt}

\end{abstract}

\maketitle

\section*{Introduction}

\medskip

Let $\La$ be an artin algebra, and denote by ${\rm mod}\hskip 0.4pt\La$ the category of finitely generated right $\La$-modules.
It is an important problem in the representation theory
of algebras to determine whether $\La$ has finite or infinite global dimension, and more
specifically, whether a simple $\La$-module has finite or infinite projective dimension. For instance, the derived category $D^b({\rm mod}\hskip 0.4pt \La)$ has Auslander-Reiten triangles if and only if $\La$ has finite global dimension; see \cite{Hap1, Hap2}. One approach to this problem is to consider the extension quiver of $\La$, which has vertices given by a complete set of non-isomorphic simple $\La$-modules and single arrows $S\to T$, where $S$ and $T$ are vertices such that ${\rm Ext}^1_\Lambda(S, T)$ is non-zero.
Then the {\it no loop conjecture} affirms that the extension quiver of $\La$ contains no loop if $\La$ is of finite global dimension, while
the {\it strong no loop conjecture}, which is due to Zacharia, strengthens this to state that a vertex in the extension quiver admits no loop if
it has finite projective dimension; see \cite{ARS, I}.

\medskip

The no loop conjecture was first explicitly established for artin
algebras of global dimension two; see \cite{GGZ}. For finite
dimensional elementary algebras,
as shown in \cite{I}, this can be easily derived from an earlier result of Lenzing on Hochschild homology in \cite{Le}. Lenzing's
technique was to extend the notion of the trace of endomorphisms of projective modules, defined by
Hattori and Stallings in \cite{Ha, Sta}, to endomorphisms of modules over a noetherian ring with finite global dimension,
and apply it to a particular kind of filtration for the regular module.

\medskip

In contrast, up to now,
the strong no loop conjecture has only been verified for
some special classes of algebras such as monomial algebras; see \cite{BF,I}, special bi\-serial algebras; see \cite{LM}, and algebras with at most two simple modules and radical cubed zero; see \cite{Je}. Many other partial results can be found in \cite{BuS, DK, GSZ, MP, Paq, Zac}. Most recently, Skorodumov generalized and localized Lenzing's filtration to indecomposable projective modules. This
allowed him to prove this conjecture for finite dimensional elementary algebras
of finite representation type; see \cite{De}.

\medskip

In this paper, we shall localize Lenzing's trace function to endomorphisms
of modules in ${\rm mod}\hskip 0.4pt \La$ with an $e$-bounded projective resolution, where $e$ is
an idempotent in $\La$. The key point is that every module in ${\rm mod}\hskip 0.4pt \La$ has an $e$-bounded projective resolution
if the semi-simple module supported by $e$ has finite injective dimension. 
This will enable us to solve the strong no loop conjecture for a large class of artin algebras including finite
dimensional elementary algebras, and particularly, for finite dimensional algebras over an algebraically closed field.

\section{Localized trace function and Hochschild homology}

\medskip

Throughout, $J$ will stand for the Jacobson radical of $\La$. The additive subgroup of $\La$ generated by the elements $ab-ba$ with $a, b\in \La$ is called the \emph{commutator group} of $\La$ and written as  $[\La,\La]$. One defines then the Hochschild homology group $\HH(\La)$ to be
$\La\hskip -0.8pt /\hskip -0.2pt[\La,\La]$. We shall say that $\HH(\La)$ is \emph{radical-trivial} if $J\subseteq [\La,\La]$.

\medskip

To start with, we recall the notion of the trace of an endomorphism $\va$ of a projective module $P$ in ${\rm mod}\hskip 0.4pt \La$, as defined by Hattori and Stallings in \cite{Ha, Sta}; see also \cite{I,Le}. Write $P=e_1\La\oplus \cdots \oplus
e_r\La$, where the $e_i$ are primitive idempotents in $\La$. Then
$\varphi=(a_{ij})_{r\times r}$, where
$a_{ij} \in e_i\La e_j$. The {\it trace} of $\varphi$ is defined to be
$$\tr(\varphi)= {\sum}_{i=1}^r a_{ii} + [\La, \La]\in \HH(\La).$$

\smallskip

We collect some well known properties of this trace function
in the following proposition, in which the property (2) is
the reason for defining the trace to be an element in $\HH(\La)$.

\medskip

\begin{PropHS} \label{tr} Let $P, P'$ be projective modules in ${\rm mod}\hskip 0.5pt\La$.

\begin{enumerate}[$(1)$]

\item If $\varphi, \psi\in {\rm End}_{\Lambda}(P)$, then
$\tr(\va + \psi)=\tr(\va)+\tr(\psi)$.

\vspace{2 pt}

\item If $\va: P \to P'$ and $\psi: P' \to P$ are $\La$-linear, then
$\tr(\va\psi) = \tr(\psi\va)$.

\vspace{2 pt}

\item If $\va=(\varphi_{ij})_{2\times 2}: P \oplus P' \to P \oplus
P'$, then $\tr(\va) = \tr(\va_{11}) + \tr(\va_{22})$.

\vspace{2 pt}

\item If $\psi: P\to P'$ is an isomorphism and $\va\in {\rm
End}_{\Lambda}(P)$,  then $\tr(\psi\va\psi^{-1})=\tr(\va)$.

\vspace{2 pt}

\item If $\va : \La \to \La$ is the left multiplication by $a \in
\La$, then $\tr(\va) = a + [\La,\La]$.

\end{enumerate}

\end{PropHS}

\medskip

Next, we recall Lenzing's extension of this notion to endomorphisms of modules of finite projective dimension. For $M\in {\rm mod}\hskip 0.5pt \La$, let $\mathcal{P}_M$ denote a projective resolution

\vspace{-3pt}
$$\cdots \longrightarrow  P_i\stackrel{d_i}{\longrightarrow} P_{i-1}\rightarrow  \cdots
\rightarrow  P_1 \stackrel{d_1}{\longrightarrow}  P_0
\stackrel{d_0}{\longrightarrow} M\rightarrow 0 \vspace{4pt}$$ of $M$ in ${\rm mod}\hskip 0.5pt \La$. For each
$\va\in {\rm End}_{\Lambda}\hskip -0.5pt (M)$, one can construct a commutative
diagram
$$\xymatrix{ \cdots \ar[r] & P_i\ar[r]^{d_i} \ar[d]^{\va_i}& P_{i-1}\ar[r]
\ar[d]^{\va_{i-1}}
& \cdots \ar[r]& P_1 \ar[r]^{d_1}\ar[d]^{\va_1}& P_0\ar[r]^{d_0} \ar[d]^{\va_0}& M\ar[r] \ar[d]^{\va}& 0\\
\cdots \ar[r] & P_i\ar[r]^{d_i}& P_{i-1}\ar[r] & \cdots \ar[r]& P_1
\ar[r]^{d_1} & P_0\ar[r]^{d_0} & M\ar[r] & 0} \vspace{2pt} $$ in ${\rm mod}\hskip 0.4pt\La$.
We shall call $\{\va_i\}_{i\ge 0}$ a {\it lifting} of $\va$ to
$\mathcal{P}_M$. If $M$ is of finite projective dimension, then one may
assume that $\mathcal{P}_M$ is bounded and define
the {\it trace} of $\va$ by

\vspace{-3pt}
$$ \tr(\va)={\sum}_{i=0}^{\infty}(-1)^i\,\tr(\va_i)\in \HH(\La), \vspace{3pt} $$
which is independent of the choice of $\mathcal{P}_M$
and $\{\va_i\}$; see \cite{Le}, and also \cite{I}.

\medskip

Our strategy is to localize this construction. Let $e$ be an idempotent in $\La$. Set
$$\La_{\hskip 0.4pt e} = \La/\La (1-e)\La.\vspace{0pt}$$
The canonical algebra projection $\La \to \La_{\hskip 0.4pt e}$
induces a group homomorphism
$$H_e : \HH(\La) \to \HH(\La_{\hskip 0.4pt e}). \vspace{2pt}$$ For an endomorphism $\va$ of a
projective module in ${\rm mod}\hskip 0.5pt \La$, we define its {\it $e$-trace} by
$$\vspace{3pt} \tr_e(\va)=H_e(\tr(\va))\in \HH(\La_{\hskip 0.4pt e}).$$

It is evident that this $e$-trace function has the properties (1) to (4) stated in Proposition \ref{tr}.
More importantly, we have the following result.

\smallskip

\begin{Lemma} \label{etr}

Let $e$ be an idempotent in $\Lambda$, and let $P$ be a projective module in ${\rm mod}\hskip 0.4pt\La$ whose top is
annihilated by $\hskip 0.2pt e$. If $\va\in {\rm End}_{\Lambda}(P)$,  then $\tr_e(\va)=0$.

\end{Lemma}

\noindent{\it Proof.} We may assume that $P$ is non-zero. Then $1-e=e_1+\cdots+e_r$, where the $e_i$ are pairwise
orthogonal primitive idempotents in $\La$. Let $\va\in {\rm End}_{\Lambda}(P)$. By Proposition
\ref{tr}(3), we may assume that $P$ is indecomposable. Then
$P\cong e_{\hskip -1pt s} \La$ for some $1\le s\le r$. By
Proposition \ref{tr}(4), we may assume that $P=e_{\hskip -1pt s} \La$. Then $\va$ is the left multiplication by some $a\in
e_{\hskip -1pt s}\La e_{\hskip -1pt s}$. By Proposition \ref{tr}(5),
$$\tr_e(\va)=H_e(a+[\La, \La])=\bar a+[\La_{\hskip 0.4pt e}, \La_{\hskip 0.4pt e}],$$ where $\bar a=a+\La\hskip 0.2pt (1-e)\La$.
Since $a=e_{\hskip -1pt
s}ae_{\hskip -1pt s}=(1-e)a(1-e)\in \La\hskip 0.2pt (1-e)\La$, we get $\tr_e(\va)=0$.  The proof of
the lemma is completed.

\medskip

To extend the $e$-trace function, we shall call a projective
resolution $\mathcal{P}_M$ of $M$ \emph{$e$-bounded} if $e$
annihilates the tops of all but finitely many terms in
$\mathcal{P}_M$. In this case, if $\va\in {\rm End}_{\Lambda}(M)$
with a lifting $\{\va_i\}_{i\ge 0}$ to $\mathcal{P}_M$ then, by
Lemma \ref{etr}, $\tr_e(\va_i)=0$ for all but finitely many $i$.
This allows us to define the {\it $e$-trace} of $\va$ by

\vspace{-5pt}

$$\tr_e(\va)={\sum}_{i=0}^{\infty}(-1)^i \,\tr_e(\va_i) \in \HH(\La_{\hskip 0.4pt e}).$$

\vspace{5pt}

\begin{Lemma} Let $e$ be an idempotent in $\La$. The $e$-trace is well defined for endomorphisms of
modules in ${\rm mod}\hskip 0.4pt\La$ having an $e$-bounded projective
resolution.
\end{Lemma}

\noindent{\it Proof.} Let $M$ be a module in ${\rm mod}\hskip 0.4pt\La$
having an $e$-bounded projective resolution

\vspace{-5pt}

$$
\mathcal{P}_M: \qquad \cdots \to P_i \stackrel{d_i}{\longrightarrow}
P_{i-1} \to \cdots \to P_1 \stackrel{d_1}{\longrightarrow} P_0
\stackrel{d_0}{\longrightarrow} M \to 0. \vspace{3pt} $$ Fix $\va\in {\rm
End}_{\Lambda}(M)$. We first show that $\tr_e(\va)$ is
independent of the choice of its lifting to $\mathcal{P}_M$. By
Proposition \ref{tr}(1), it amounts to proving that $\sum_{i=0}^\infty (-1)^i\,\tr_e(\va_i)=0$ for any
lifting $\{\va_i\}_{i\ge 0}$  of the zero endomorphism of $M$. Indeed, let $h_i : P_i \to P_{i+1}$ be
morphisms such that $\va_0 = d_1h_0$
and $\va_i= d_{i+1}h_i + h_{i-1}d_i$, for $i \ge 1.$ Applying
Proposition \ref{tr}, we get

\vspace{-8pt}

$$\tr_e(\va_i)
= \tr_e(d_{i+1}h_i) + \tr_e(h_{i-1}d_i)= \tr_e(d_{i+1}h_i) +
\tr_e(d_ih_{i-1}),$$

\vspace{2pt}

\noindent for $i\ge 1$. On the other hand, by assumption, there exists some $m \ge
0$ such that $e$ annihilates the top of $P_i$ for $i \ge m$. By
Lemma \ref{etr}, $\tr_e(d_{m+1}h_{m})=0$ and $\tr_e(\va_i)=0$ for $i
\ge m$. This yields

\vspace{-8pt}

\begin{eqnarray*}
  {\sum}_{i=0}^{\infty}(-1)^i \,\tr_e(\va_i) &=& \tr_e(\va_0) + {\sum}_{i=1}^{m}(-1)^i \,\tr_e(\va_i) \\
  &=& \tr_e(d_1h_0) + {\sum}_{i=1}^{m}(-1)^i\left(\tr_e(d_{i+1}h_i) +
  \tr_e(d_ih_{i-1})\right) \\
 &=& (-1)^{m}\,\tr_e(d_{m+1}h_{m})\\
  &=&0.
\end{eqnarray*}

Next, we show that $\tr_e(\va)$ is independent of the choice of the $e$-bounded projective
resolution $\mathcal{P}_M$. Suppose that $M$ has another $e$-bounded
projective resolution
$$\mathcal{P}'_M:
\qquad \cdots \to P'_i \stackrel{d'_i}{\longrightarrow} P'_{i-1} \to
\cdots \to P'_1 \stackrel{d'_1}{\longrightarrow} P'_0
\stackrel{d'_0}{\longrightarrow} M \to 0.$$ 
Considering $\va$, we get morphisms $u_i : P_i \to P'_i$ with $i\ge 0$ such that $d'_0u_0= \va d_0$ and $d'_iu_i=u_{i-1}d_i$ for $i \ge 1$. Similarly, considering $1_M$, we obtain maps
$v_i: P_i'\to P_i$ with $i\ge 0$ such that $d_0v_0= d'_0$ and $d_iv_i=v_{i-1}d'_i$ for
$i \ge 1$. Observe that $\{v_iu_i\}_{i\ge 0}$ and $\{u_iv_i\}_{i\ge 0}$ are liftings of $\va$
to $\mathcal{P}_M$ and $\mathcal{P}_M'$, respectively. By Proposition \ref{tr}(2),
we have $$\qquad {\sum}_{i=0}^{\infty}(-1)^i\,\tr_e(u_iv_i)=
{\sum}_{i=0}^{\infty}(-1)^i \,\tr_e(v_iu_i).$$  The proof of
the lemma is completed.

\bigskip

In the sequel, $S_e$ will stand for
the semi-simple $\La$-module $e\La/eJ$. Suppose that $S_e$ has finite injective
dimension. If $M$ is a module in ${\rm mod}\hskip
0.3pt\La$, then $\Ext_{\La}^{i}(M,S_e)=0$ for all sufficient large integers $i$, that is, the minimal projective resolution of $M$ is $e$-bounded. Therefore, the $e$-trace is defined for every endomorphism in ${\rm mod}\hskip 0.4pt\La$. In particular, if $\La$ is of finite global dimension, then we recover Lenzing's
trace function by taking $e=1_A$.

\medskip

\begin{Prop}\label{additivity} Let $e$ be an idempotent in $\La$.
Consider a commutative diagram
$$\xymatrixrowsep{18pt}
\xymatrix{0 \ar[r] & L\ar[r]^u\ar[d]^{\va_{_L}} & M \ar[r]^v\ar[d]^{\va_{_M}} & N\ar[d]^{\va_{_N}} \ar[r] & 0\\
0 \ar[r] & L \ar[r]^u & M \ar[r]^v & N \ar[r] &0}$$ in ${\rm mod}\hskip 0.4pt\La$ with exact rows. If $L, N$ have $e$-bounded projective
resolutions, then so does $M$ and $\tr_e(\va_{_M}) = \tr_e(\va_{_L}) +
\tr_e(\va_{_N}).$
\end{Prop}
\smallskip
\noindent{\it Proof.} Assume that $L$ and $N$ have $e$-bounded projective resolutions as follows$\,:$

\vspace{-5pt}

$$\mathcal{P}_L: \qquad \cdots \to P_i \stackrel{d_i}{\longrightarrow} P_{i-1}
\to \cdots \to P_1 \stackrel{d_1}{\longrightarrow} P_0
\stackrel{d_0}{\longrightarrow} L \to 0$$ and
$$\mathcal{P}_N:\qquad  \cdots \to P'_i
\stackrel{d'_i}{\longrightarrow} P'_{i-1} \to \cdots \to P'_1
\stackrel{d'_1}{\longrightarrow} P'_0
\stackrel{d'_0}{\longrightarrow} N \to 0. \vspace{4pt}$$ By the
Horseshoe lemma, there exists in ${\rm mod}\hskip
0.5pt \La$ a commutative diagram

\vspace{-8pt}
$$\xymatrixrowsep{18pt} \xymatrixcolsep{18pt}\xymatrix{
 \cdots \ar[r]& P_i\ar[r]^-{d_i} \ar[d]^{q_i}& P_{i-1} \ar[r] \ar[d]^{q_{i-1}}&
\cdots \ar[r] & 
P_0 \ar[r]^{d_0}\ar[d]^{q_{_0}}& L \ar[r]\ar[d]^u& 0\\
 \cdots \ar[r]& P_i\oplus P_i'\ar[r]^-{d''_i} \ar[d]^{p_i}& P_{i-1} \oplus P'_{i-1}\ar[r] \ar[d]^-{p_{i-1}}&
\cdots \ar[r] & 
P_0\oplus P'_0
\ar[r]^-{d''_0}\ar[d]^{p_{_0}}& M \ar[r]\ar[d]^v& 0\\
 \cdots \ar[r]& P'_i\ar[r]^-{d'_i} & P'_{i-1} \ar[r] &
\cdots \ar[r] & 
P'_0 \ar[r]^-{d'_0}& N \ar[r]& 0\vspace{5pt}}$$  with exact rows, where $q_i={1 \choose {\hskip 0.4pt 0\hskip 0.5pt}}$, $p_i=(0,1)$
for all $i\ge 0$. In particular, the middle row
is an $e$-bounded projective resolution of $M$ which we denote by
${\mathcal P}_M$. Choose a lifting $\{f_i\}_{i\ge 0}$ of $\va_{_L}$
to $\mathcal{P}_L$ and a lifting $\{g_i\}_{i\ge 0}$ of $\va_{_N}$ to
$\mathcal{P}_N$. It is well known; see, for example, \cite[p.
46]{We} that there exists a lifting $\{h_i\}_{i\ge 0}$ of $\va_{_M}$
to $\mathcal{P}_M$ such that
$$\xymatrixrowsep{18pt}
\xymatrix{0 \ar[r] & P_i \ar[r]^-{q_i} \ar[d]^{f_i} & P_i \oplus
P'_i \ar[r]^-{p_i}\ar[d]^{h_i} & P'_i \ar[r] \ar[d]^{g_i} & 0 \\
0 \ar[r] & P_i \ar[r]^-{q_i} & P_i \oplus P'_i \ar[r]^-{p_i} &
P'_i \ar[r] & 0 } \vspace{3pt} $$ is commutative, for every $i\ge 0$. Since
$h_iq_i=q_if_i$ and $g_ip_i=p_ih_i$, we can write $h_i$ as a $(2
\times 2)$-matrix whose diagonal entries are $f_i$ and $g_i$. Thus
$\tr_e(h_i) = \tr_e(f_i) + \tr_e(g_i)$ by Proposition \ref{tr}(3).
As a consequence, $\tr_e(\va_M) = \tr_e(\va_N) + \tr_e(\va_L).$
The proof of the proposition is completed.

\bigskip

Finally, we shall describe the Hochschild homology group
$\HH(\La_{\hskip 0.4pt e})$ in case $S_e$ has finite injective dimension.

\medskip

\begin{Theo} \label{maintheo} Let $\La$ be an artin algebra, and let $e$ be an idempotent in $\La$.
If $S_e$ has finite injective dimension, then $\HH(\La_{\hskip 0.4pt e})$ is
radical-trivial.
\end{Theo}

\noindent{\it Proof.} Suppose that $S_e$ has finite injective
dimension. Then the $e$-trace is defined for every endomorphism in
${\rm mod}\hskip 0.5pt \La$. Let $x \in \La$ be such that $\bar
x=x+\La(1-e)\hskip 0.3pt\La \,$ lies in the radical of $\La_{\hskip
0.4pt e}$, which is $(J + \La(1-e)\La) /\La(1-e)\La$. Hence, $\bar x = \bar a$ for some
$a\in J$. Let $r>0$ be such that $a^r=0$, and consider the chain
$$0=M_r\subseteq M_{r-1}\subseteq \cdots \subseteq M_1\subseteq M_0=\La,$$
of submodules of $\La$, where $M_i = a^i\La$, $i=0, \ldots, r$. Let
$\va_0: \La\to \La$ be the left multiplication by $a$. Since
$\va_0(M_i)\subseteq M_{i+1}$, we see that $\va_0$ induces
morphisms $\va_i : M_i \to M_i$, $i=1, \ldots, r$, such that
$$\xymatrixrowsep{17pt}
\xymatrix{0 \ar[r] & M_{i+1}\ar[r]\ar[d]^{\va_{i+1}} & M_i \ar[r]\ar[d]^{\va_i} & M_i/M_{i+1}\ar[d]^{0} \ar[r] & 0\\
0 \ar[r] & M_{i+1} \ar[r] & M_i \ar[r] & M_i/M_{i+1} \ar[r] &0}$$
commutes, and hence $\tr_e(\va_i) = \tr_e(\va_{i+1})$ by Proposition \ref{additivity}, for $i=0, 1, \ldots, r-1$. Applying
Proposition \ref{tr}(5), we get
$$\bar a+[\La_{\hskip 0.4pt e}, \La_{\hskip 0.4pt e}]=H_e(a+[\La, \La])=H_e({\rm tr}(\va_0))=\tr_e(\va_0) =\tr_e(\va_r)=0,$$
that is, $\bar x = \bar a\in [\La_{\hskip 0.4pt e},\La_{\hskip 0.4pt
e}]$.  The proof of the theorem is completed.

\medskip

Taking $e=1_A$, we recover the
following well known result; see, for example, \cite{Le}.

\medskip

\begin{Cor}
If $\La$ is an artin algebra of finite global dimension, then
$\HH(\La)$ is radical-trivial.
\end{Cor}

\medskip

Indeed, if $\La$ is a finite dimensional algebra of finite global dimension over
a field of characteristic zero, then all the Hochschild homology groups ${\rm H\hspace{-1pt}H}_i(\La)$ with
$i \ge 1$ vanish; see \cite{Le}. However, in the situation as in Theorem \ref{maintheo}, the higher Hochschild homology
groups of $\La_{\hskip 0.2pt e}$ do not necessarily vanish and $\La_{\hskip 0.4pt e}$ may be of infinite global
dimension.

\bigskip

\noindent {\sc Example.}
Let $\La = kQ/I$, where $k$ is a field, $Q$ is the quiver
$$\xymatrixrowsep{19pt}\xymatrixcolsep{19pt}\xymatrix{1 \ar[r]^{\alpha}\ar[d]_{\gamma} & 2 \ar[d]^{\beta}\\ 4 \ar[r]^{\delta}
& 3 \ar[ul]_{\varepsilon} }$$ and $I$ is the ideal in $kQ$
generated by $\alpha \beta - \gamma \delta, \beta \varepsilon,
\delta \varepsilon, \varepsilon\alpha$. One can show that $\La$ has
finite global dimension. Now, let $e$ be the sum of the primitive
idempotents in $\La$ corresponding to the vertices $1,2,3$. Then
$\La_{\hskip 0.4pt e}$ is a Nakayama algebra with radical squared
zero, which clearly has infinite global dimension. By Theorem \ref{maintheo}, $\HH(\La_{\hskip 0.4pt e})$ is
radical-trivial. However, a direct computation shows that ${\rm
H\hspace{-1pt}H}_2(\La_{\hskip 0.4pt e})$ is non-zero; see also
\cite{IZ}.

\section{Main results}

\medskip

The main objective of this section is to apply the previously obtained result to solve the strong no loop conjecture for finite dimensional algebras over an algeb\-raically closed field. We start with an artin algebra $\La$ with a primitive idempotent $e$. We shall say that $\La$ is {\it locally commutative} at $e$ if $e\La e$ is commutative and that $\La$ is {\it locally
commutative} if it is locally commutative at every primitive
idempotent. Moreover,  $e$ is called {\it basic} if $e\La$ is not isomorphic to any direct summand of
$(1-e)\La$. In this terminology, $\La$ is basic if and only if all its primitive idempotents are
basic.

\medskip

\begin{Theo}\label{main}
Let $\La$ be an artin algebra, and let $e$ be a basic primitive
idempotent in $\La$ such that $\La/J^2$ is locally commutative
at $e+J^2$. If $S_e$ has finite projective or injective dimension,
then $\Ext^1_{\Lambda}(S_e,S_e)=0$.
\end{Theo}

\noindent{\it Proof.} Firstly, we assume that $S_e$ is of finite injective dimension. For proving that
$\Ext^1_{\Lambda}(S_e,S_e)=0$, it suffices to show that $e J e
/eJ^2e = 0$. Let $a \in eJe$. Then $a +
\La(1-e)\La \in [\La_{\hskip 0.4pt e}, \La_{\hskip 0.4pt e}]$ by Theorem \ref{maintheo}. Since
$e$ is basic, $e\La(1-e)\La e \subseteq eJ^2e$.  This yields an
algebra homomorphism $$f: \La_e  \to e\La e/eJ^2e: x+\La(1-e)\La \mapsto exe + eJ^2e.$$ Thus, $a +
eJ^2e=f(a + \La(1-e)\La)$ lies in the commutator group of $e\La e/eJ^2e$. On the other hand, $e\La e/eJ^2e\cong
(e+J^2)(\La/J^2)(e+J^2)$, which is commutative. Therefore, $a +
eJ^2e=0$, that is, $a \in eJ^2e$. The result follows in this case.

Next, assume that $S_e$ has finite projective dimension. Let $D$
be the standard duality between $\mmod \La$ and $\mmod \hskip 0.3pt
\La^{\, \rm op}$. Then $D(S_e)$ is the simple $\La^{\,\rm
op}$-module supported by the idempotent $e^{\rm o}$ corresponding to $e$, which is of finite injective
dimension. Observe that the quotient of $\La^{\rm op}$ modulo its
radical square is also locally commutative at the class of $e^{\rm o}$ modulo the radical square. By what
we have proven,
$\Ext_{\La}^1(S_e,S_e)\cong \Ext_{\La^{\, \rm op}}^1(D(S_e),D(S_e))=0.$
The proof of the theorem is completed.

\bigskip

\noindent{\sc Remark.} The preceding result establishes the strong no loop conjecture for
basic artin algebras $\La$ such that $\La/J^2$ is locally
commutative.

\bigskip

Now we shall specialize this result to finite dimensional algebras over a field.
Recall that such an algebra is called {\it elementary} if its simple modules are all one dimensional over the base field; see \cite{ARS}.

\medskip

\begin{Theo}\label{main2}

Let $\La$ be a finite dimensional algebra over a field $k$, and let
$S$ be a simple $\La$-module which is one dimensional over $k$. If $S$ has finite
projective or injective dimension, then $\Ext^1_{\Lambda}(S,S)=0$.

\end{Theo}

\noindent{\it Proof.} Let $e\in \La$ be the primitive idempotent supporting $S$. Then $\La$
has a complete set $\{e_1, \ldots, e_n\}$ of primitive orthogonal idempotents with $e=e_1$. We may assume that $e_1\La, \ldots, e_r\La$, with $1\le r\le n$, are the non-isomorphic indecomposable projective modules in ${\rm mod}\hskip 0.5pt \La$. Then
$$\La/J \cong M_{n_1}(D_1) \times \cdots \times M_{n_r}(D_r),$$
where $D_i={\rm End}_{\Lambda}(e_i\La/e_iJ)$ and $n_i$ is the number of indices $j$ with $1\le j\le n$ such that $e_j\hskip 0.4pt \La\cong e_i\La$, for $i=1, \ldots, r.$ Now $S$ is a simple $M_{n_1}(D_1)$-module, and hence $S\cong D_1^{n_1}$. Since $S$ is one dimensional over $k$, it is one dimensional over $D_1$. In particular, $n_1=1$. That is, $e$ is a basic primitive idempotent. Moreover,
$e\La e/eJe \cong Se \cong k$. Thus, for $x_1,x_2\in e \La e$, we can write
$x_i=\lambda_ie + a_i$, where $\lambda_i\in k$ and $a_i\in eJe$,
$i=1, 2$. This yields $x_1x_2-x_2x_1=a_1a_2-a_2a_1\in eJ^2e.$ Therefore,
$e\La e/ eJ^2e$ is commutative, and so is $(e+J^2)(\La/J^2)(e+J^2)$.
The result follows immediately from Theorem \ref{main}. The proof of the
theorem is completed.

\medskip\vspace{3pt}

\noindent{\sc Remark.} The preceding theorem establishes the strong no loop conjecture for finite dimensional elementary algebras, and hence for finite dimensional algebras over an algebraically closed field.

\medskip\vspace{3pt}

We shall extend our results in this direction. Let $\La$ be a finite dimensional elementary algebra over a field $k$. We may assume that $\La=kQ/I$, where $Q$ is a finite quiver, $kQ$ is the path algebra, and $I$ is an admissible
ideal in $kQ$; see \cite{ARS}. Recall that $I$ is {\it admissible} if $(kQ^+)^n\subseteq I\subseteq (kQ^+)^2$ for some $n\ge 2$, where $kQ^+$ is the ideal in $kQ$ generated by the arrows, and {\it monomial} if $I$ in addition is generated by some paths. In this setting, the extension quiver of $\La$ is isomorphic to the quiver obtained from $Q$ by shrinking the possible multiple arrows. If $p_1, \ldots, p_r$ are distinct paths in $Q$ of length $\ge 2$ from one vertex to another, then a $k$-linear combination
$$\rho=\lambda_1p_1+\cdots+\lambda_rp_r\vspace{3pt} $$
is called a \emph{minimal
relation} for $\La$ if $\rho\in I$ and $\sum_{i \in J}\lambda_ip_i \not \in I$ for any  $J\subset \{1,\cdots,r\}$.  Moreover, let $\sigma = \alpha_1 \alpha_2 \cdots \alpha_r$ be an
oriented cycle in $Q$, where the $\alpha_i$ are arrows. The support
of $\sigma$, written as ${\rm supp}(\sigma)$, is the set of vertices in
$Q$ occurring as starting points of $\alpha_1, \ldots, \alpha_r$. The \emph{idempotent supporting $\sigma$} is the sum of
all primitive idempotents in $\La$ associated to the vertices in
${\rm supp}(\sigma)$. Write
$$\sigma_1 = \sigma, \; \sigma_i = \alpha_{i} \cdots \alpha_r\alpha_1 \cdots \alpha_{i-1}, \; i=2, \ldots, r, $$
called the {\it cyclic permutations} of $\sigma$.
We shall say that $\sigma$ is {\it cyclically free} in $\La$ if none of the $\sigma_i$ with $1\le i\le r$ is a summand of a minimal relation for $\La$, and {\it cyclically non-zero} in $\La$ if  none of the $\sigma_i$ lies in $I$.

\medskip

\begin{Theo} \label{main2} Let $\La = kQ/I$ with $Q$ a finite quiver and $I$ an admissible ideal in $kQ$, and let
$\sigma$ be an oriented cycle in $Q$ with supporting idempotent $e\in \La$. If $\sigma$ is cyclically free in $\La$, then $S_e$ has
infinite projective and injective dimensions.

\end{Theo}

\noindent{\it Proof.} Suppose that $\sigma$ is cyclically free in $\La$.  If $\sigma$ is a power of a shorter oriented cycle $\delta$, then it is easy to see that $\delta$ is also cyclically free in $\La$ and ${\rm supp}(\delta)={\rm supp} (\sigma)$. Hence, we may assume that $\sigma$ is not a power of any shorter oriented cycle. Let $\sigma_1, \ldots, \sigma_r$, where $\sigma_1=\sigma$, be the cyclic permutations of $\sigma$. It is then well known that the $\sigma_i$ with $1\le i\le r$ are pairwise distinct.

For any $p\in kQ$, denote by $\tilde p$ its class in
$\La$ and by $\bar{p}$ the class of $\tilde p$ in $\La_{\hskip 0.4pt e}$.
Let $W$ be the vector subspace of $\La_{\hskip 0.4pt e}$ spanned by the classes
$\bar p$, where $p$ ranges over the paths in $Q$ different from
$\sigma_1, \ldots, \sigma_r$. Then, there exist paths $p_1, \ldots,
p_m$ in $Q$ such that $\{\bar p_1, \ldots, \bar p_m\}$ is a
$k$-basis of $W$. We claim that $\{\bar \sigma_1, \ldots, \bar
\sigma_r, \bar p_1, \ldots, \bar p_m\}$ is a $k$-basis of
$\La_{\hskip 0.4pt e}$. Indeed, it clearly spans $\La_{\hskip
0.4pt e}$. Assume that
$${\sum}_{i=1}^r\lambda_i \bar \sigma_i + {\sum}_{j=1}^m\nu_j\bar
p_j=\bar 0, \; \lambda_i, \nu_j\in k.$$  That is, $\sum
\lambda_i \tilde \sigma_i + \sum \nu_j \tilde p_j \in \La (1-e)
\La$.  Then
$${\sum}_{i=1}^r\lambda_i \,\tilde \sigma_i +
{\sum}_{j=1}^m\nu_j\,\tilde p_j = {\sum}_{l=1}^s \mu_l \,\tilde
q_l, \;\mu_l\in k,$$ where $q_1, \ldots, q_s$ are distinct paths in $Q$ passing through a vertex not in ${\rm supp}(\sigma)$. Fix some
$t$ with $1\le t\le r$. Letting $\varepsilon_t$ be the trivial path in $Q$
associated to the starting point $a_t$ of $\sigma_t$, we get
$${\sum}_{i=1}^r \lambda_i \varepsilon_t\,\sigma_i\,\varepsilon_t + {\sum}_{j=1}^m \nu_j \,\varepsilon_t\, p_j \,\varepsilon_t- {\sum}_{l=1}^s \mu_l \varepsilon_t\, q_l \,\varepsilon_t \in I.$$

Note that the non-zero elements of
the $\varepsilon_t{\sigma}_i\varepsilon_t$,  $\varepsilon_t
p_j \varepsilon_t$, $\varepsilon_t \, q_l
\,\varepsilon_t\in kQ$ are distinct oriented cycles from
$a_t$ to $a_t$. Since $\sigma$ is cyclically free in $\La$, we have $\lambda_j=0$ whenever $\varepsilon_t\sigma_j\varepsilon_t$ is non-zero. In particular,
$\lambda_t=0$. Therefore, the $\lambda_i$ are all
zero, and so are the $\nu_j$. This proves our claim. Suppose now that
$\bar \sigma  \in [\La_{\hskip 0.4pt e},\La_{\hskip 0.4pt e}]$. Then
\begin{equation} \bar{\sigma} = {\sum}_{i=1}^n\eta_i(\bar{u}_i\bar{v}_i - \bar{v}_i\bar{u}_i)
\end{equation}
where $\eta_i \in k$ and $u_i, v_i\in \{\sigma_1, \ldots,
\sigma_r,  p_1, \ldots,  p_m\}$. For each $1\le i\le n$, we see easily that
$u_iv_i\not\in \{\sigma_1, \ldots, \sigma_r\}$ if and only if
$v_iu_i\not\in \{\sigma_1, \ldots, \sigma_r\}$, and in this case, $\bar{u}_i\bar{v}_i- \bar v_i\bar
u_i\in W$. Therefore, the equation $(1)$ becomes
\begin{equation} \bar{\sigma} = {\sum} \,\eta_{ij}(\bar{\sigma}_i - \bar{\sigma}_j) + w,
\end{equation} where $\eta_{ij}\in k$ and $w\in W$.
Let $L$ be the linear form on $\La_{\hskip 0.4pt e}$, which sends
each of $\bar{\sigma}_1, \ldots, \bar \sigma_r$ to $1$ and vanishes on $W$. Since $\sigma=\sigma_1$, applying $L$ to
the equation $(2)$ yields $1 = 0$, a contradiction. Therefore, the class of $\bar \sigma$ in $\HH(\La_{\hskip 0.4 pt e})$ is non-zero. Since $\bar \sigma$ lies in the radical of $\La_{\hskip 0.4pt e}$, by Theorem \ref{maintheo}, $S_e$ has infinite projective and injective dimensions.
The proof of the theorem is completed.

\medskip \vspace{5pt}

\noindent{\sc Example.} Let $\La
=kQ/I$, where  $Q$ is the following quiver
$$\xymatrixrowsep{25pt}\xymatrixcolsep{25pt}
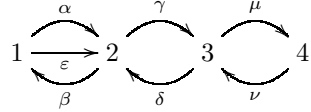
\xymatrix{1 \ar@/^1pc/[r]^(0.5){\alpha} \ar[r]_{\varepsilon} & 2
\ar@/^1pc/[r]^(0.5){\gamma} \ar@/^1pc/[l]^(0.5){\beta} & 3
\ar@/^1pc/[l]^{\delta} \ar@/^1pc/[r]^{\mu} & 4 \ar@/^1pc/[l]^{\nu}
} \vspace{-3pt} $$ and $I$ is the ideal in $kQ$ generated by
$\alpha\beta,\delta\gamma,\beta\varepsilon,\varepsilon\beta,
\nu\delta,\nu\mu,\mu\nu,\gamma\mu$,
$\alpha\gamma\delta\beta\alpha\gamma-\varepsilon\gamma$. It is easy to see that the oriented cycle $\beta\alpha \gamma\delta$ is cyclically free in $\La$. By
Theorem \ref{main2}, one of the simple modules $S_1, S_2, S_3$ has infinite projective
dimension.

\medskip

\begin{Cor} \label{cor1} Let $\La = kQ/I$ with $Q$ a finite quiver and $I$ an admissible ideal in $kQ$. If $Q$ contains an oriented cycle which is cyclically free in $\La$, then $\La$ has infinite global dimension.
\end{Cor}

\medskip

If $I$ is a monomial ideal in $kQ$, then an oriented cycle in
$Q$ is cyclically free in $\La$ if and only if it is cyclically non-zero in $\La$.
This yields the following consequence, which can also be derived from results in \cite{IZ}.

\medskip

\begin{Cor} Let $\,\La=kQ/I$ with $Q$ a finite quiver and $I$ a monomial ideal in $kQ$. If $Q$ contains an oriented cycle which is cyclically non-zero in $\La$, then $\La$ has infinite global dimension.
\end{Cor}

\medskip

To conclude, we would like to draw the reader's attention to an even stronger version
of the no loop conjecture as follows.

\medskip

\begin{Conj}  Let $S$ be a simple module over an artin algebra. If $\Ext^1(S,S)$ is
non-zero, then $\Ext^i(S,S)$ is non-zero for infinitely many integers $i$.
\end{Conj}

\medskip

This conjecture was originally posed under the name of {\it extreme
no loop conjecture} in \cite{LM}. It remains open except for monomial
algebras and special biserial algebras; see \cite{GSZ, LM}.

\medskip
\bigskip

{\sc Acknowledgment.} This research was supported in part by the
Natural Scien\-ces and Engineering Research Council of Canada and was
carried out during the first author's visit to the Universit\'e de Sherbrooke.

\end{document}